\documentclass[11pt,letterpaper]{amsart}

\usepackage{cmap}  
\usepackage[T2A]{fontenc}
\usepackage[utf8]{inputenc}
\usepackage[english]{babel}
\usepackage[usenames]{color}
\usepackage{amsmath,amssymb,amsthm}
\usepackage{arcs}
\usepackage{enumitem}
\usepackage{epsfig}
\usepackage{filecontents}
\usepackage{graphicx}
\usepackage[pdftex,colorlinks=true,linkcolor=blue,urlcolor=red,unicode=true,hyperfootnotes=false,bookmarksnumbered]{hyperref}
\usepackage{ifthen}
\usepackage{import}
\usepackage{indentfirst}
\usepackage{mathtools}
\usepackage{cancel}
\usepackage{xcolor}

\usepackage[margin=1in, total={8.5in, 11in}]{geometry}

\newcommand{\eps}{\varepsilon}

\renewcommand{\le}{\leqslant}
\renewcommand{\ge}{\geqslant}

\newcommand{\ff}{\mathcal{F}}
\newcommand{\s}{\mathcal{S}}
\newcommand{\g}{\mathcal{G}}

\newcommand{\aaa}{\mathcal{A}}

\newcommand{\h}{\mathcal H}
\newcommand{\uu}{\mathcal U}
\newcommand{\m}{\mathcal }

\newcommand{\T}{\mathcal{T}}
\newcommand{\W}{\mathcal{W}}

\newtheorem{theorem}{Theorem}

\newtheorem{lemma}[theorem]{Lemma}
\newtheorem{corollary}[theorem]{Corollary}

\newtheorem{observation}[theorem]{Observation}

\title{A complete $t$-intersection theorem for families of spanning trees} 
\author{Elizaveta Iarovikova}
\address{Moscow Institute of Physics and Technology, Russia; Email: {\tt liza.fm@yandex.ru}}
\author{Andrey Kupavskii}
\address{Moscow Institute of Physics and Technology, St. Petersburg State University, Innopolis University, Russia; Email: {\tt kupavskii@ya.ru}}
\date{}

\begin{document}
\maketitle
\begin{abstract}
    Let $\m T_n$ denote the set of all labelled spanning trees of $K_n$. A family $\ff \subset \T_n$ is \textit{$t$-intersecting} if for all $A, B \in \ff$ the  trees $A$ and $B$ share at least $t$ edges. In this paper, we determine for $n>n_0$  the size of the largest $t$-intersecting family $\ff\subset T_n$ for all meaningful values of $t$ ($t\le n-1$). This result is a rare instance  when a complete $t$-intersection theorem for a given type of structures is known. 
\end{abstract}
\section{Introduction}
We use standard notation $[n]$ for the set $\{1,\ldots, n\}$, $2^X$ for the power set of $X$ and ${X\choose k}$ for the set of all $k$-element subsets of $X$. A family of sets $\ff$ is {\it $t$-intersecting} if for any $A,B\in \ff$ we have $|A\cap B|\ge t.$ For shorthand, we say that a family is {\it intersecting} instead of $1$-intersecting. The study of $t$-intersecting families goes back to the seminal paper of Erd\H os, Ko and Rado \cite{EKR}, who studied intersecting and $t$-intersecting families in $2^{[n]}$ and ${[n]\choose k}$. Shortly after, Katona \cite{Kat} determined the largest $t$-intersecting family of sets in $2^{[n]}$ for all $n,t.$ In particular, he showed that for even $n+t$, the largest such family consists of all sets of size at least $(n+t)/2$.  

Determining the largest $t$-intersecting family for ${[n]\choose k}$ proved to be a much harder task. After several partial results \cite{F1, W, FF3} by Frankl, Wilson, and Frankl and F\"uredi, Ahlswede and Khachatrian \cite{AK} confirmed the conjecture of Frankl \cite{F1} and showed that, for any $n,k,t$, the largest $t$-intersecting family in ${[n]\choose k}$ must, up to isomorphism, have the form 
\begin{equation}\label{exfra} \Big\{F\in {[n]\choose k}: |F\cap[t+2i]|\ge t+i\Big\}  \end{equation}
for some integer $i\ge 0$. A similar result was shown by Filmus \cite{Film} for $2^{[n]}$ with the $p$-biased measure.

Over the last 65 years, the Erd\H os--Ko--Rado type theorems were obtained for many different structures. However, the only other case where a complete $t$-intersection theorem is known that we are aware of is the result of Frankl and Wilson \cite{FW86} that determines the largest family of $k$-dimensional subspaces of an $n$-dimensional vector space over a finite field $GF(q)$ with $q$ elements. There, we say that the family of subspaces is $t$-intersecting if the intersection of any two subspaces from the family has dimension at least $t$.

One of the most notable and well-studied structures in that respect is $t$-intersecting permutations. We say that two permutations $\pi,\sigma$ intersect in at least $t$ elements if there are $t$ distinct elements $x_1,\ldots, x_t$ such that $\sigma(x_i) = \pi(x_i)$ for all $t$. Frankl  and Deza showed that any intersecting family of permutations $[n]\to[n]$ has size at most $(n-1)!$. After a series of results by Ellis, Friedgut and Pilpel \cite{EFP}, Kupavskii and Zakharov \cite{spreads}, Keller, Lifshitz, Minzer and Sheinfield \cite{KLMS}, Kupavskii \cite{permutations} determined the largest family of $t$-intersecting permutations for all $t\le (1-\epsilon)n$ and $n>n_0(\epsilon)$. For $t\ge n/2$, the correct answer is no longer $(n-t)!,$ corresponding to the family of all permutations that fix the elements in $[t],$ but rather the examples of the type \eqref{exfra}.

Permutations $[n]\to [n]$ can be thought of as perfect matching in the complete bipartite graph $K_{n,n}.$ Two permutations $t$-intersect iff the corresponding matchings share at least $t$ edges. In this paper, we study $t$-intersecting families of closely related objects: spanning trees. We show that the situation with $t$-intersection theorems for trees is somewhat simpler than that for permutations.
Let $\T_n$ denote the set of labelled spanning trees in a complete graph on $n$ vertices. Cayley's tree formula states that $|\T_n| = n^{n-2}$. We say that a family $\ff \subset \T_n$ is \textit{$t$-intersecting} if for all $A, B \in \ff$ the  trees $A$ and $B$ share at least $t$ edges. 

A common example of $t$-intersecting families are \textit{trivial} $t$-intersecting families. Those are the families of sets that all contain a given $t$-element set. In this paper, we show that such trivial intersecting families are always extremal for $t$-intersecting families of trees.
The main result of this paper is the following theorem.

\begin{theorem}\label{thmbigt}
   There exists an absolute constant $n_0$ such that for all $n \ge n_0$ and $2 \le t \le n-2$ the following holds. If $\ff\subset \T_n$ is $t$-intersecting then $|\ff| \le c_{n, t}n^{n-2-t}$, where $c_{n, t}$ is the largest possible product of $n-t$ positive integers with sum $n$. Moreover, the equality holds if and only if $\ff$ is the family of all spanning trees containing a fixed forest $F$ with $t$ edges and connected components of almost equal sizes.
\end{theorem}

To prove these results, we use the spread approximation technique introduced in \cite{spreads} and developed in subsequent papers of the second author and coauthors \cite{Kup55, permutations, FKhr, Fedya}. In particular, our proof builds on the proof of the main result  in \cite{permutations} by the second author. 

Note that each tree has $n-1$ edges, and thus an $(n-1)$-intersecting family of spanning trees can contain at most one spanning tree. Thus the bound $t\le n-2.$ Our result extends the result of Frankl, Hurlbert, Ihringer, Kupavskii, Lindzey, Meagher and Tej Pantangi, who showed the same result for  $n \ge 2^{19}$ and $1 < t \leqslant \frac{n}{4032\log_2n}$ this result was proved in \cite{EKR_trees}. In this paper, for clarity we only present the proof for  large $t$, even though in \cite{EKR_trees} the authors also used spread approximations and so the proof could be unified. The authors made certain conjectures concerning the size of the largest $t$-intersecting families of spanning trees for large $t,$ and our main result  confirms these conjectures.

The case $t = 1$ is somewhat different, as the largest family is not exactly trivial. The following was proved in  \cite{EKR_trees}. 

\begin{theorem}\label{thmt1}
  For $n > 2^{19}$, let $\ff$ be $1$-intersecting family of spanning trees on $n$ vertices. Then $$|\ff| \le 2n^{n-3} + n - 2.$$ Moreover, equality holds if and only if $\ff$ consists of all stars plus all trees containing a fixed edge.
\end{theorem}

We should note that our proof of Theorem~\ref{thmbigt} implies that for $t<(1-\epsilon)\frac n2$ or $t\ge n/2$ any family that has size $(1+o(1))$ times the size of the extremal family must be very close to a trivial $t$-intersecting family. In the regime $(1-\epsilon)\frac n2\le t<n/2,$ however, there is another example of approximately the same size: the analogue of the family as in \eqref{exfra} with $i=1.$ For permutations, such example overtakes the trivial family once $t$ passes across $n/2.$ Here, however, once $t$ becomes larger than $n/2,$ the size of the corresponding $i=1$ example becomes constant times smaller than the size of the trivial example, and the ratio grows as $t$ grows. This could serve as one explanation for the fact that we managed to prove the $t$-intersection theorem in the entire range.

In the next section, we discuss spread families. In Section~\ref{sec3} we compare the sizes of $t$-intersecting families of trees in the spirit of \eqref{exfra} and show that trivial $t$-intersecting families are always the largest. The following sections are dedicated to the proof of the main theorem.
\section{Preliminaries}
We treat spanning trees as $(n-1)$-element subsets of the set of all edges, therefore $\T_n \subset {{n\choose 2}\choose n-1}$. We use the following notation that is standard when working with families of sets. Given families $\ff$, $\s$ and a set $X$, we define
$$\ff(X) :=\{ A\setminus X  :\    A \in \ff, \,  X \subset A \}$$
$$\ff[X] :=\{ A  : \ A\in \ff, \,  X \subset A \}$$
$$\ff[ \s ] := \bigcup\limits_{X \in \s }\ff[X].$$

Given $q, r >1$ a family $\ff\subset2^{[n]}$ is said to be \textit{$r$-spread} if $$|\ff(X)| < r^{-|X|}|\mathcal{F}|$$ for all $X \subset [n]$.

We will use several observations concerning spread families.

\begin{observation}\label{obs_spread_1}
  If for some $\alpha>1$ we have a family $\ff\subset {[n]\choose \le k}$ such that $|\ff| > \alpha^k$, then there is a set $X$ such that $|\ff(X)|>1$  and $\ff(X)$ is $\alpha$-spread.
\end{observation}

\begin{proof}
    If $\ff$ is $\alpha$-spread then we are done. Otherwise, find an inclusion-maximal set $X$ such that $|\ff(X)| \ge \alpha^{-|X|}|\ff|.$ Note that in this case $|\ff(X)| \ge \alpha^{-k}|\ff| > 1$. Assume $\ff(X)$ is not spread, then there is a nonempty set $Y$ disjoint from $X$ such that $|\ff(X \cup Y)| \ge \alpha^{-|Y|}|\ff(X)| \ge \alpha^{-|X|-|Y|}|\ff|$. But then we should have chosen $Y \cup X$ in place of $X$ which was supposed to be inclusion-maximal.
\end{proof}
 Arguing as in the proof of Observation \ref{obs_spread_1}, we get the following observation.
\begin{observation}\label{obs_spread_2}
  Fix $\alpha>1$ and consider a family $\ff$. If $X$ is an inclusion-maximal set with the property that $\ff(X) \ge \alpha^{-|X|}|\ff|$, then $\ff(X)$ is $\alpha$-spread.
\end{observation}

The main theorem concerning spread families is the following statement from \cite{Alw} called the `Spread lemma'. We will need its sharpening due to Hu \cite{Hu}. Later, Stoeckl in \cite{Sto} established the same theorem  with $1 + h_2(\delta) \le 2$ instead of constant $2$, where $h_2(\delta) = -\delta \log_2 \delta -(1 - \delta) \log_2(1 - \delta)$ is the binary entropy of $\delta$.

\begin{theorem}\label{spread_lemma}
  If for some $n,k,r\ge 1$ a family $\ff\subset {[n]\choose \le k}$ is $r$-spread and $W$ is an $(\beta\delta)$-random subset of $[n]$, then $$\Pr[\exists F\in \ff\ :\ F\subset W]\ge 1-\Big(\frac 2{\log_2(r\delta)} \Big)^\beta k.$$
\end{theorem}

We will also use an application of this theorem to find disjoint sets in spread families.

\begin{lemma}\label{find_disjoint}
Let $\g_1, \g_2 \subset 2^{[m]}$ be $R$-spread uniform families of uniformity at most $k$. Suppose that $R>2^{5}\log_22k$. Then there exist disjoint sets $G_1 \in \g_1, G_2 \in \g_2$.
\end{lemma}

\begin{proof}
Let us put $\beta= \log_2(2k)$ and $\delta = (2\log_2(2k))^{-1}$. Note that $\beta\delta = \frac 12$ and $R\delta > 2^{4}$ by our choice of $R$. Theorem~\ref{spread_lemma} implies that a $\frac{1}2$-random subset $W_i$ of $[m]$ contains a set from $\g_i$ with probability strictly bigger than
$$1-\Big(\frac 2{\log_2 2^{4}}\Big)^{\log_2 (2k)} k = 1-2^{-\log_2 (2k)} k = \frac 12.$$
Note that $[n]\setminus W$ is also a $\frac12$-random subset of $[m]$. A simple union bound implies that with positive probability there are $G_1 \in \g_1$ and $G_2 \in \g_2$ such that $G_1 \subset W, G_2 \subset [m]\setminus W$. Thus there are disjoint $G_1 \in \g_1, G_2 \in \g_2$.
\end{proof}

\section{Reconstruction of spanning trees}\label{sec3}
In this section, we discuss families of spanning trees containing a fixed subforest. The main goal being to compare the size of different examples of $t$-intersecting families and ultimately show that trivial families are always the largest. 

We will use a generalization of Cayley's tree formula from~\cite[Lemma 6]{tree_counter} that counts the number of labelled trees that contain a specific subforest. We will also compare the sizes of various $t$-intersecting families. We refer to the number of edges in the forest as to its \textit{size}.

\begin{lemma}\label{count_trees}
Let $F$ be a spanning forest in the $K_n$ that has $m$ connected components of sizes $q_1, q_2, \dots, q_m$. Then the number of spanning trees in $K_n$, that contain $F$ is exactly
\[
q_1 q_2 \cdots q_m \ n^{n-2-\sum_{i=1}^m (q_i-1)} .
\]
\end{lemma}

Let $c_{n, t}$ denote the maximal size of $|\T_n[F]|/n^{n-2-t}$ over all spanning forests of complete $n$-vertex graph with $t$ edges, that is $$c_{n, t} = \max_{\sum_{i=1}^m(q_i -1) = t, \sum_{i=1}^{m}q_i = n}q_1\cdot\ldots\cdot q_m.$$ First, note that $m = n-t$. It is easy to see that the maximum value of $c_{n, t}$ is attained if all variables $q_i$ are almost equal, that is, differ by no more than 1.

More precisely, if $r$ is such non-negative integer that $n = (n-t)\lfloor\frac{n}{n-t}\rfloor + r, r < n-t$. Then we have \begin{equation}\label{eqc} c_{n, t} = \Big\lceil\frac{n}{n-t}\Big\rceil^r\cdot\Big\lfloor\frac{n}{n-t}\Big\rfloor^{n-t-r}.\end{equation}
Recall that the elements of $\{q_i\}_{i=1}^m$ correspond to the sizes of connected components of a certain spanning forest $F$. Therefore a forest $F$ maximizes the value of $\T_n[F]$ among all forests with $t$ edges if and only if it consists of $n-t$ connected components of sizes $\lfloor\frac{n}{n-t}\rfloor$ and $\lceil\frac{n}{n-t}\rceil$. We note that the structure inside the components plays no role here.

Later we will prove that the extremal example of a $t$-intersecting family is always a trivial family, i. e. a family that consists of all spanning trees containing a fixed spanning forest $F$. Note that in the extremal case $F$ should be as described above. However we will start with a broader collection of $t$-intersecting families. For $r = 0, 1, \ldots, n - 1 - t$ and a spanning forest $F$ of size $t+2r$ consider a $t$-intersecting family $$\uu_{n, t, r, F} = \{A\in \T_n \mid |A\cap F|\ge t+r\}.$$ 

Denote $$\aaa_{n, t, r} = \max_{F -\text{forest of size } t+2r} |\uu_{n, t, r, F}|.$$ (We remark that here, as an exception, we use a calligraphic letter to denote a quantity rather than a family.) Note that $$\aaa_{n, t, 0} = c_{n, t}n^{n - 2 - t}.$$ One of the main goals of this section is to show that with fixed $n, t$ the value of $\aaa_{n, t, 0}$ is the largest among $\{\aaa_{n, t, r}\}_{r = 0}^{n-1-t}$. This result will be proved in Lemma \ref{antr_relative} and Lemma \ref{ant1}.

We start with comparing $c_{n, t+r}$ and $c_{n, t}$, that is we prove the following lemma, which we then use as a base  $r = 1$ case of induction on $r$.
    
\begin{lemma}\label{cnt_relative_base}
    Let $n,t$ be nonnegative integers such that $t\le n - 2$. Then 
    \begin{enumerate}[label=(\roman*)]
        \item $c_{n, t+1}/c_{n, t} \le 2$;
        \item if $t \ge 2n/3$ then $c_{n, t+1}/c_{n, t} \le e/3$.
        \item if $n/2 \le t < 2n/3$ then $c_{n, t+1}/c_{n, t} \le 9/8$.
        \item $c_{n, t+1}/c_{n, t} \ge 1/n$
    \end{enumerate}
\end{lemma}
\begin{proof}
     Recall that $c_{n,t}$ can be thought of as of the product of $n-t$ integer factors of almost equal values  with sum $n$, cf. \eqref{eqc}. Consider a sequence of such factors $q_1 \le q_2 \le \ldots \le q_{n-t}$. For $i = 1, \ldots, n - t - 1$ put $p_i = q_{i+1}$. Then for $j = 1, \ldots, q_1$ increase $p_j$ by one where the numeration of $p_i$'s is considered to be cyclic, $p_{n-t} = p_1$. It means that if $q_1 > n - t -1$, then there will be some $p_i$'s that we increase by $1$ more then once.

    It is easy to verify that the $p_1, \ldots, p_{n-t-1}$ obtained during this procedure are $n-t-1$ integers that differ by at most $1$ and that have sum $n$. Thus $c_{n, t+1} = p_1\cdot \ldots \cdot p_{n-t-1}$. Put $k = \lfloor \frac{n}{n-t}  \rfloor = q_1$. Before this process \begin{equation}\label{eq_pq}p_1\cdot\ldots\cdot p_{n-t-1}/q_1\cdot\ldots\cdot q_{n-t} = 1/q_1 = 1/k.\end{equation} Each change $p_i := p_i + 1$ increases the numerator on the LHS of the equation \eqref{eq_pq} by the factor of at most $1 + 1/k$. There are exactly $k$ such steps, therefore

    $$\frac{c_{n, t+1}}{c_{n, t}} = \frac{p_1\cdot\ldots\cdot p_{n-t-1}}{q_1\cdot\ldots\cdot q_{n - t}} \le \frac{\left(1 + \frac1k\right)^k}{k} = f(k).$$ 

    Estimate (i) in the statement of the lemma can be verified directly for $k=1, 2$., i.e. we should check that $f(1) \le 2$, $f(2) \le 2$. For $k \ge 3$ we have $f(k) < e/3$ since $$\left(1 + \frac1k\right)^k < e < 3\leqslant k,$$ thus $f(k) \le e/3$. To see (ii) note that if $t \ge 2n/3$ then $n-t < n/3$ and $k = \lfloor \frac n{n-t}\rfloor \ge 3$. Lastly $n/2 \le t < 2n/3$ implies $k=2$ and $f(2) = 9/8$. Note that since $f(1) = 2$ and $f(2) = 9/8$ presented bounds for $t \le 2n/3$ are best possible.

    Lastly, since at start of the procedure $p_i = q_{i + 1}$ and  during the procedure $p_i$'s can only increase, we get $$\frac{c_{n, t+1}}{c_{n, t}} = \frac{p_1\cdot\ldots\cdot p_{n-t-1}}{q_1\cdot\ldots\cdot q_{n - t}} \ge \frac1{q_1} \ge \frac1n,$$ which implies (iv) and completes the proof.
\end{proof}

Having the base case $r=1$, we proceed by induction and prove the following statement.

\begin{lemma}\label{cnt_relative}
    Let $n,t, r$ be nonnegative integers such that $t + r\le n - 2$. Then 
    \begin{enumerate}[label=(\roman*)]
        \item $c_{n, t+r}/c_{n, t} \le 2^r$;
        \item if $t \ge 2n/3$ then $c_{n, t+r}/c_{n, t} \le (e/3)^r$.
        \item if $n/2 \le t < 2n/3$ then $c_{n, t+r}/c_{n, t} \le (9/8)^r$.
        \item $c_{n, t+r}/c_{n, t} \ge 1/n^r$
    \end{enumerate}
\end{lemma}

\begin{proof}
    Case $r = 0$ is obvious, case $r = 1$ was treated in Lemma \ref{cnt_relative_base}. Now we proceed by induction on~$r$. Assume the statement is true for $r = 1, \ldots, s$. Consider $r = s+1$.

    To prove (i) we use (i) of Lemma \ref{cnt_relative_base}: $$\frac{c_{n, t+s+1}}{c_{n, t}} = \frac{c_{n, t+s+1}}{c_{n, t+s}} \cdot \frac{c_{n, t+s}}{c_{n, t}} \le 2\cdot 2^{s} = 2^{s+1}.$$

    The proof of (ii) and (iii) is similar. If $t \ge 2n/3$, then $t + s + 1 \ge 2n/3$ and $$\frac{c_{n, t+s+1}}{c_{n, t}} = \frac{c_{n, t+s+1}}{c_{n, t+s}} \cdot \frac{c_{n, t+s}}{c_{n, t}} \le (e/3)\cdot (e/3)^{s} = (e/3)^{s+1}.$$ If $n/2 \le t < 2n/3$ then $t+s \ge n/2$. Due to (ii) and (iii) of Lemma \ref{cnt_relative_base} we have $c_{n, t+s+1}/c_{n, t+s} \le 9/8$. Thus $c_{n, t+s+1}/c_{n, t} \le (9/8)^{s+1}$.

    Lastly, using the same argument and (iv) of Lemma \ref{cnt_relative_base}, we get (iv) of the statement.
\end{proof}

Using this lemma, we can compare $\aaa_{n, t, r}$ and $\aaa_{n, t, 0}$. 
\begin{lemma}\label{antr_relative}
    There is $n_0 \in \mathbb{N}$ and an absolute constant $\alpha < 1$ such that for any positive integers $n,t, r$ such that $t\le n-2$, $t+r\le n-1$ and $n\ge n_0$ we have $\aaa_{n, t, r} \le \alpha \aaa_{n, t, 0}$ whenever we have either $r \ge 2$, or $r=1$ and $t\ge n/2$.
\end{lemma}

\begin{proof}
    For any $n, t, r$ and forest $F$ with $n$ vertices and $t+2r$ edges the size of $\uu_{n, t, r, F}$ is not bigger than ${t+2r\choose t+r}\aaa_{n, t+r, 0}$ by a simple union bound. Taking the maximum over all such $F$ we get $$\aaa_{n, t, r} \le {t+2r\choose t+r}\aaa_{n, t+r, 0} = {t+2r\choose r}c_{n, t+r}n^{n-2-t-r} \leqslant 2^rc_{n, t}{t+2r\choose r}n^{n-2-t-r},$$ where the last inequality is due to (i) of Lemma \ref{cnt_relative}.

    There is a positive integer $r_0$ such that for all integer $r \ge r_0$ we have $4^r/r! \le  4^{r_0}/r_0! < 1$. Put $\alpha_1 = 4^{r_0}/r_0!$. From the inequality above and the fact that for all $r \ge r_0$ we have $t+2r \le 2n$, we get that the following holds: \begin{equation}\label{eq_antr}\aaa_{n, t, r} \le 2^rc_{n, t}\cdot\frac{(2n)^r}{r!}\cdot n^{n-2-t-r} \leqslant \frac{4^r}{r!}\cdot c_{n, t}n^{n-2-t} \leqslant \alpha_1\aaa_{n, t, 0}.\end{equation}

    Now consider $r < r_0$. Here we treat cases $t < n/2, n/2\le t<2n/3$ and $t \ge 2n/3$ separately. 
    
    \textbf{Case I: $t < n/2$.} This is the only case where we exclude $r=1$.  Take $\eps > 0$ such that $\alpha_2 = (1+\eps)^2/2 < 1$. Then $(1+\eps)^r/r! < 1$ holds for all $r \ge 2$. There is such $n_1$ that for all $n > n_1$ it is true that $t+2r \le n/2 + 2r_0 < (1 + \eps)n/2.$ We write a bound similar to the inequality \eqref{eq_antr}. $$\aaa_{n, t, r} \le 2^rc_{n, t}\cdot\frac{((1+\eps)n/2)^r}{r!}\cdot n^{n-2-t-r} \leqslant \frac{(1+\eps)^r}{r!}\cdot c_{n, t}n^{n-2-t} \leqslant \alpha_2\aaa_{n, t, 0}.$$

    \textbf{Case II: $n/2 \le t < 2n/3$.} In this case we choose $\eps > 0$ so that $\alpha_3 = (1+\eps)\cdot 3/4 < 1$. There is $n_2$ such that for any $n > n_2$ holds $t+2r < (1+\eps)2n/3$. We also use (iii) of Lemma \ref{cnt_relative}. 

    $$\aaa_{n, t, r} \le (9/8)^rc_{n, t}\cdot\frac{((1+\eps)\cdot 2/3\cdot n)^r}{r!}\cdot n^{n-2-t-r} \leqslant \frac{(1+\eps)^r\cdot(3/4)^r}{r!}\cdot c_{n, t}n^{n-2-t} \leqslant \alpha_3\aaa_{n, t, 0}.$$
    
    \textbf{Case III: $2n/3\le t$.} In the this case we use (ii) of Lemma \ref{cnt_relative}. Choose $\eps > 0$ such that $\alpha_4 = (1 + \eps)e/3 < 1$. For $n > n_3$, where $n_3$ is large enough, $t+2r < n + 2r_0 < (1 + \eps)n$. Using the same bound as in the two previous cases, 
    $$\aaa_{n, t, r} \le (e/3)^r c_{n, t}\cdot\frac{(1+\eps)^rn^r}{r!}\cdot n^{n-2-t-r} \le \frac{((1+\eps)e/3)^r}{r!}\aaa_{n, t, 0} < \alpha_4\aaa_{n, t ,0}.$$

    Therefore the statement of the lemma is true for $0 < \alpha = \min\{\alpha_1, \alpha_2, \alpha_3, \alpha_4\} < 1, n_0 = \max\{n_1, n_2, n_3\}$. Note that we needed $r \ge 2$ only in the first case.
    
\end{proof}

In case $t < n/2$ we do not have such a good bound for $\aaa_{n, t, 1}$, however we can still compare it to $\aaa_{n, t, 0}$.

\begin{lemma}\label{ant1}
    There is a positive constant $n_0$ such that for any positive integers $n,t$ such that $t<n/2$ and $n > n_0$ we have $\aaa_{n, t, 0} > \aaa_{n, t,1}$.
\end{lemma}

\begin{proof}
    In this case $\aaa_{n, t, 0} = 2^tn^{n-2-t}$. Note that for some $A, |A| = t+2$ the value of $\aaa_{n, t, 1}$ can be bounded by a simple bound which counts some sets more than once.
    \begin{align}
    \label{eq_antr1}\aaa_{n, t, 1} &\leqslant\sum_{B \subset A, |B| = t+1}|\T_n[B]| \\
     \label{eq_antr2}&\le  (t+2)\cdot c_{n, t+1} n^{n-2-t-1} \\
     \label{eq_antr3}&\le 2^{t+1}(t+2)n^{n-2-t-1} \leqslant \aaa_{n, t, 0} \cdot \frac {2(t+2)}{n}.
    \end{align}

    This expression implies $\aaa_{n, t, 1} < \aaa_{n, t, 0}$ for $t < n/2 - 2$. To complete the proof we should treat four remaining cases: $$t \in \{n/2 - 2, n/2 - 3/2, n/2 - 1, n/2 - 1/2\}.$$
    
    In case $t = n/2 - 2$ by the series of inequaitites \eqref{eq_antr1} we get $\aaa_{n, t, 1} \le \aaa_{n, t, 0}$ with equality if and only if all $B's$ are matchings. It implies that $A$ does not contain a cycle, so $\T_n[A]$ is non-empty. In  the RHS of \eqref{eq_antr1} we count all trees from $\T_n[A]$ $t+2$ times, and thus \eqref{eq_antr1} is strict. This guarantees $\aaa_{n, t, 1} < \aaa_{n, t, 0}$.
    
    In the three remaining cases, we have $|A| = t+2 > n/2$, so at least two edges in $A$ have a common vertex. Only two $(t+1)$-element subsets of $A$ do not contain one of these two edges. Therefore at most two $(t+1)$-element subsets of $A$ are matchings. All other $(t+1)$-element subsets have at least one connected component of size $3$ or greater. Thus, the size of the corresponding $\m T_n[B]$ is at most  $3\cdot 2^{t-1}n^{n-2-t-1}<2^{t+1}n^{n-2-t-1}=c_{n,t+1}n^{n-2-t-1}$.  It implies that the inequality \eqref{eq_antr2} is strict. 
\end{proof}

The previous two lemmas together show that $\aaa_{n, t, 0}$ is the biggest among $\{\aaa_{n, t, r}\}_{r = 0}^{n-1-t}$.

We will also use a lemma from \cite{EKR_trees} that bounds the number of spanning trees containing a fixed forest $F$ and avoiding a spanning tree $T$ outside of $F$.

\begin{lemma}\label{avoiding}
    Let $F$ be a forest on $n$ vertices with $t$ edges, $t \le n - 55$, $T$ be a spanning tree on $n$ vertices which is not a star, $F \not\subset T$. Then the number of spanning trees in $\T_n$ that contain $F$ and avoid edges of $T$ outside of $F$ is at least $n^{n-t-27}$.
\end{lemma}

We prove separately the same result when $F$ is a star. Here we are interested in trees that intersect a fixed star in less than $t$ edges.

\begin{lemma}\label{avoiding_star}
    Let $F$ be a forest on $n$ vertices with $t$ edges, $1 < t \le n - 3$, $S$ be a complete star on $n$ vertices, $F \not\subset S$. Then the number of spanning trees in $\T_n$ that contain $F$ and have less than $t$ common edges with $S$ is greater than $(n-1)^{n-2-t}$.
\end{lemma}
\begin{proof}
    Let $v$ be such a vertex in $S$ that is connected with all other vertices.

    \textbf{Case I: $v$ is isolated in $F$.} There are at least $(n-1)^{(n-1) - 2 - t}$ trees that cover all vertices in $[n]\setminus \{v\}$ and contain $F$. For each such tree, add an edge between $v$ and any other vertex. Thus we get a tree that has exactly $1 < t$ common edge with $S$. During this procedure, we constructed at least $n(n-1)^{n-3-t} > (n-1)^{n-3-t}$ with properties as described in the statement.

    \textbf{Case II: $v$ is not isolated in $F$.} There are at least $(n-1)^{(n-1) - 2 - (t-1)}$ trees that cover all $n$ vertices but the vertex $v$ and contain all edges of $F$ except those incident to $v$. If we add all edges from $F$ to such a tree, we get a tree that has less then $t$ common edges with $S$. Therefore, the lemma is proved in both cases.
\end{proof}

\section{Outline of the proof}
Given the result of \cite{EKR_trees}, we may assume that $t > 2n^{1 - \eps/7}$ for a very small $\eps$ to be defined later. Let $\ff$ be the largest $t$-intersecting family of spanning trees of a complete $n$-vertex graph. Assume $|\ff| \ge \aaa_{n, t, 0} = c_{n, t}n^{n-2-t}$. We will prove that in this case $|\ff| = \aaa_{n, t, 0}$  and there is such a forest $F$ with $t$ edges that $\ff = \T_n[F]$. Moreover, $F$ is a union of $n-t$ connected components of sizes  $\lfloor\frac{n}{n-t}\rfloor$ and $\lceil\frac{n}{n-t}\rceil$. 

We cannot immediately use the spread approximation method in order to find a $t$-intersecting approximation of low uniformity. First, for some $t'<t$ we find a $t'$-intersecting set $\s$ of forests of sizes smaller than $(1 + \eps)t$ such that it covers most of the family $\ff$. Then we use the simplification-peeling technique to analyse the structure of $\s$ and find a restriction $X$ such that $\ff(X)$ is `dense' in $\T_n(X)$. Inside this $\ff(X)$ we identify a very spread subfamily, and then via a greedy procedure we build a $t$-intersecting spread approximation $\s'$.

In the fourth step, we apply a variant of the peeling-simplification process to $\s'$ in order to determine the most important 'layer'. As a result, we prove that all but $o(|\ff|)$ sets of $\ff$ lie in some $\uu_{n, t, r, F}$. The last step is to show that in the extremal family the remainder must be empty. To this end, we argue indirectly and identify a spread subfamily inside  the remainder $\ff\setminus \uu_{n, t , r, F}$. This spread subfamily  forbids too many sets from $\uu_{n, t , r, F}$, implying that if $\ff$ is extremal then  $\ff = \uu_{n, t, r, F}.$

Let us comment separately on the case of very large $t$ (say, $t>0.999n$. We treat this case together with the rest, but, actually, if spelled out separately, it would actually be simpler. We do not need to do any kind of spread approximation and may go directly to the simplification-peeling procedure. 

\section{First spread approximation}

We begin the proof of the main result by constructing the first, `rough' spread approximation. In the next stage, we shall repeatedly apply this theorem to $\ff$ and some of its `dense' subfamilies.

\begin{theorem}\label{1SA}
  There is such $\eps_0 > 0$ that for any $\eps < \eps_0, n > n_0(\eps)$ the following holds. Let $t\ge1$ be an integer such that $n^{1-\eps/3} \le t \le n- 2$. Consider a family $\ff\subset \T_n$ that is $t$-intersecting.
  Put $t'= t-n^{1-\epsilon/3}$, $q' = \min\{(1 + \eps)t, n-1\}$. Then there exists a  $t'$-intersecting family $\mathcal S\subset \T_n^{\le q'}$ of spanning forests of size at most $q'$ and a family $\ff'$ such that the following holds.
  \begin{itemize}
    \item[(i)] $\ff\setminus \ff'\subset \T_n[\s]$;
    \item[(ii)] for any $B\in \s$ there is a family $\ff_B\subset \ff$ such that $\ff_B(B)$ is $n^{\epsilon/2}$-spread;
    \item[(iii)] $|\ff'|\le n^{-\frac 13\eps t}\cdot c_{n, t}n^{n-t-2}$.
  \end{itemize}
\end{theorem}

\begin{proof}[Proof of Theorem~\ref{1SA}] We construct $\mathcal S$ as follows. Put $r:=
n^{\eps/2}$ and consider the following procedure for $i=1,\ldots $ with $\ff^1:=\ff$.
\begin{enumerate}
    \item Find a maximal $S_i$ that  $|\ff^i(S_i)|\ge  r^{-|S_i|}|\ff^i|$.
    \item If $|S_i|> (1+\eps)t$ or $\ff^i = \emptyset$ then stop. Otherwise, put $\ff^{i+1}:=\ff^i\setminus \ff^i[S_i]$.
\end{enumerate}
Note that Observation~\ref{obs_spread_2} and maximality of $S_i$ imply that $\ff^i(S_i)$ is $r$-spread.
Let $N$ be the step of the procedure for $\ff$ at which we stop. The family $\s$ is defined as follows: $\s:=\{S_1,\ldots, S_{N-1}\}$. Clearly, $|S_i|\le q'$ for each $i\in [N-1]$. The family $\ff_{B}$ promised in (ii) is defined to be $\ff^i[S_i]$ for $B=S_i$. Next, note that if $\ff^N$ is non-empty, then $(1+\eps)t < |S| \le n-1$. Put $a = |S| - t > \eps t$. Then
\begin{align*}
|\ff^N|\le&r^{|S_N|}  |\ff^{N}(S_N)|\le  r^{|S_N|}c_{n, |S_n|}n^{n-2-|S_n|} \le \\ 
\overset{\ref{cnt_relative}, (i)}{\le}&
n^{\frac 12\eps|S_N|}\cdot 2^a\cdot c_{n, t}\cdot n^{n - 2 - t - a} = \\
=& c_{n, t}n^{n-2-t}n^{\eps t/2 + \eps a/2 - a+a\log_n2} \le \\
\le& n^{-\frac13\eps t}\cdot c_{n, t} \cdot n^{n-2-t}. 
\end{align*}
The last inequality is because $$a\left(\frac \eps2 - 1 + \log_n2\right) + \frac{\eps t}2 \leqslant -\frac{\eps t}3 \Longleftrightarrow$$ $$\Longleftrightarrow a\ge \frac{5 \eps t}{6(1 - \eps/2 - \log_n2)}.$$ which holds for $\eps < \eps_0$ and $n > n_0(\eps)$ if we fix $\eps_0$ small enough and $n_0(\eps)$ large enough, because $$a \ge \eps t \ge  \frac{5 \eps t}{6(1 - \eps/2 - \log_n2)}.$$ We are only left to verify that $\s$ is $t'$-intersecting if we chose $\eps_0$ and $n_0(\eps)$ appropriately. We prove it below.

\begin{lemma}\label{lemtint_1SA} The family $\mathcal S$ is $t'$-intersecting. \end{lemma}
\begin{proof}

  Take two (not necessarily distinct) $A_1,A_2\in \mathcal S$  and assume that $|A_1\cap A_2|<t'$. Recall that the families $\ff_{A_i}(A_i)$, $i=1,2$ are both $r$-spread with $r = n^{\eps/2}$. 
  For $i\in [2]$ put $j:=3-i$ and put
  $$\g_i:=\ff_{A_i}(A_i)\setminus\Big\{F\in \ff_{A_i}(A_i): |F\cap A_j\setminus A_i|\ge \frac{n^{1-\eps/3}}2\Big\}.$$ The size of the latter family is at most
 \begin{align*}{|A_j\setminus A_i|\choose n^{1-\epsilon/3}/2}\max_{X: |X| = n^{1-\eps/3}/2, X\cap A_i = \emptyset}|\ff_{A_i}(A_i\cup X)|\le&\\  {n\choose n^{1-\eps/3}/2}n^{-\frac\eps 4 n^{1-\eps/3}}|\ff_{A_i}(A_i)|\le&\\
  (2e n^{\eps/3})^{\frac 12n^{1-\eps/3}}n^{-\frac\eps 4 n^{1-\eps/3}}|\ff_{A_i}(A_i)|\le&\\
  n^{\frac{\eps}5n^{1 - \eps/3} - \frac{\eps}4n^{1 - \eps/3}}|\ff_{A_i}(A_i)| \le&\frac 12 |\ff_{A_i}(A_i)|.\end{align*}
  Two last inequalities require $n$ to be large enough, the exact minimal value of $n$ depends on $\eps$ only. We can guarantee that by the choice of $n_0(\eps)$. This implies that $|\g_i|\ge \frac 12 |\ff_{A_i}(A_i)|$.
 Because of this and the trivial inclusion $\g_i(Y)\subset \ff_{A_i}(A_i\cup Y)$, valid for any $Y$,  we conclude that $\g_i$ is $ \frac r2$-spread for both $i\in [2]$. Since $n$ can be sufficiently large depending only on $\eps$, we have $\frac r2 > 2^{11}\log_2(2n)$.

  What follows is an application of Theorem~\ref{spread_lemma}. Let us put $\beta= \log_2(2n)$ and $\delta = (2\log_2(2n))^{-1}$. Note that $\beta\delta = \frac 12$ and $\frac r2\delta > 2^{4}$ by our choice of $r$.  Theorem~\ref{spread_lemma} implies that a $\frac{1}2$-random subset $W_i$ of ${{[n]\choose 2}\choose 2}\setminus A_i$ contains a set from $\g_i$ with probability strictly bigger than
  $$1-\Big(\frac 2{\log_2 2^{4}}\Big)^{\log_2 (2n)} n = 1-2^{-\log_2 (2n)} n = \frac 12.$$

  Consider a random partition of ${[n]\choose 2}\setminus (A_1\cup A_2)$ into $2$ parts $U_1',U_2'$, where $\Pr[x\in U_i']=1/2$ independently for each $i\in [2]$ and $x\in {[n]\choose 2}\setminus (A_1\cup A_2)$. Next, for $i\in[2]$ and $j=3-i$, put $U_i:= U_i'\cup (A_j\setminus A_i)$. Then  $U_i$ is a random subset of the same ground set as $W_i$ above, moreover, the probability of containing each element is at least that for $W_i$. It implies that there is $F_i \in \g_i$ such that $F_i\subset U_i$ with probability strictly bigger than $\frac 12$. Using the union bound, we conclude that, with positive probability, it holds that there are such $F_i$, $F_i\subset U_i,$ for both  $i \in[2]$. Fix such choices of $U_i$ and $F_i$, $i \in [2]$. Then, on the one hand, each $F_i\cup A_i$ belongs to $\ff$ and, on the other hand, $|(F_1\cup A_1)\cap (F_2\cup A_2)| =|A_1\cap A_2|+|F_1\cap A_2|+|F_2\cap A_1|< |A_1\cap A_2|+n^{1-\epsilon/2}<t$. The first inequality follows from the definition of $\g_i$, and the second inequality due to $|S_1\cap S_2|<t'$ and the definition of $t'$. This is a contradiction with $\ff$ being $t$-intersecting.
  \end{proof}
This concludes the proof of Theorem~\ref{1SA}.\end{proof}

\section{Peeling procedure}\label{sec_peel}
In this section, we analyse the `peeling' procedure that captures the structure of the $t$-intersecting families of spanning forests. One of the key objectives is to control the number of trees coming from higher uniformity spanning forests.

 We say that a $t$-intersecting family $\h$ is {\it maximal} if for any $H \in \h$ and any proper subset $X \subsetneq H$ there exist $H' \in \h$ such that $|H\cap H'|<t$ and, moreover, $H_1\not\subset H_2$ for $H_1,H_2\in \h$ (i.e., $\h$ is an {\it antichain}).

\begin{observation}\label{obs_peeling}
For any positive integers $n,q$ and a $t$-intersecting family $\s \subset \T^{\le q}_n$, there exists a maximal $t$-intersecting family $\h$  such that for every $H\in \h$ there is $S\in \s$ such that $H\subseteq S$. Consequently, for any family $\ff\subset \T_n$ we have $\ff[\s] \subset \ff[\h]$.
\end{observation}
The easiest way to see it is to construct $\h$ by repeatedly replacing sets in $\s$ by their proper subsets while preserving the $t$-intersecting property.

Consider a $t$-intersecting family $\s\subset \T^{\le q}_n$. Let us iteratively define the following series of families.
\begin{enumerate}
    \item Put $\h_{q-t}=\s$.
    \item For $k = q-t,q-t-1, \ldots, 0$ we put $\W_k = \h_k \cap {[n] \choose t+k}$ and let $\h_{k-1}$ be the family given by Observation~\ref{obs_peeling} when applied to $\h_{k}\setminus \W_{k}$ playing the role of $\mathcal S$.
\end{enumerate}
Thus, we `peel' the sets of the largest uniformity in $\h_k$ and replace the resulting family by a maximal $t$-intersecting family using Observation~\ref{obs_peeling}.  Remark that $\h_k$ is $t$-intersecting for each $k=q-t,q-t-1\ldots,0$ by definition. We summarize the basic properties of these series of families in the following lemma.

\begin{lemma}\label{lemma_peel} The following properties hold for each $k = q-t,q-t-1\ldots, 0$
\begin{itemize}
  \item[(i)] All sets  in $\h_k$ have size at most $t+k$.
  \item[(ii)] For $k\le q-t-1$ and any $\aaa \subset \T_n$ we have $\aaa[\h_{k+1}]\subset \aaa[\h_{k}]\cup \aaa[\W_{k+1}]$.
  \item[(iii)] There is no set $X$ of size $\le t+k-1$
      such that $\W_k(X)$ (or any subfamily of $\h_k(X)$) is $\alpha$-spread for $\alpha>k + 1$. 
\end{itemize}
\end{lemma}
\begin{proof}
(i) This easily follows by reverse induction on $k$. The base case is that all sets in $\s$ have size at most $q$, and the induction step is via the definition of $\h_{k-1}$.

(ii) We have $\aaa[\h_{k+1}] = \aaa[\h_{k+1} \setminus \W_{k+1}] \cup \aaa[\W_{k+1}]$ and, by the definition of $\h_{k}$, we have $\aaa[\h_k]\supset \aaa[\h_{k+1}\setminus \W_{k+1}]$.

(iii) Assume that there is such a set $X$. By the maximality of $\h_k$, there must be a set $F\in \h_k$ such that $|X\cap F|=z<t$.
At the same time, the family $\h_k$ is  $t$-intersecting. So, for each set  $G\in \W_k[X]$ we have $|G\cap F|\ge t$, and thus $|F\cap (G\setminus X)| \ge t-z$ for each such $G$. Thus, $\W_k[X]\subset \cup_{S\in {F\choose t-z}}\W_k[X\cup S]$ and so, using the $\alpha$-spreadness of $\W_k(X)$, we get \begin{align*}|\W_k(X)|\le&\ {|F|-z\choose t-z}\max_{F'\in {F\choose t-z}}|\W_k(X\cup F')|\\
\le&\ {|F|-z\choose t-z}\alpha^{-(t-z)}|\W_k(X)| \\
=&\ \prod_{i=1}^{t-z}\frac{|F|-t+i}i\alpha^{-(t-z)}|\W_k(X)|\\
\le&\ (|F|-t+1)^{t-z}\alpha^{-(t-z)}|\W_k(X)|<|\W_k(X)|,\end{align*}
where the last inequality is due to $|F|-t+1\le k + 1<\alpha$. This a contradiction.
\end{proof}

Our next goal is to bound the size of $\W_k$.
Consider two sets $A,B\in \h_k$ such that $A\cap B = I$, $|I|=t$. (There are such sets by maximality of $\h_k$.) Recall that the sizes of sets of $A,B$ are at most $t+k$.
Any set $C\in \W_k$ intersects both $A$ and $B$ in at least $t$ elements. Denote $C_0=C\cap I$, $C_1=C\cap (A\setminus I)$ and $C_2 = C\cap (B\setminus I)$. Then, putting $|C_0|=t-j$, we must have $|C_1|, |C_2|\ge j$. Select $C_1'\subset C_1,C_2'\subset C_2$ such that $|C_1'| = |C_2'| = j$. The family $\W_k(C_0\cup C_1'\cup C_2'\cup X)$ by Lemma~\ref{lemma_peel} (iii) is not  $\alpha$-spread subfamily for $\alpha>k + 1$ and any $X$, and thus $|\W_i(C_0\cup C_1'\cup C_2')|\le (k+1)^{k-j}$ by Observation~\ref{obs_spread_1}. From here, we see that we can upper bound the number of sets in $\W_k$ as follows:
$$|\W_k|\le \sum_{j=0}^{k}{t\choose t-j}{|A|-t\choose j}{|B|-t\choose j} (k+1)^{k-j}\le \sum_{j=0}^{k}f(j),$$
where
\begin{equation}\label{eqfj} f(j):= {t\choose t-j}{k\choose j}^2 (k+1)^{k-j}.\end{equation}

The observation from \cite{permutations} below is based on simple calculations that check whether the${f(j)/}{f(j+1)}$ is bigger than $1$ or not.
\begin{observation}\label{w_bounds} For sufficiently large $t$, $t\ge k$, we have the following. Denote by $j_0$ the value at which $f(j)$ attains its maximum.
  \begin{enumerate}
    \item For $k$ such that $t=o(k^2)$, we have $j_0 =(1+o(1))(tk)^{1/3}$.
    \item For $k=\Theta (\sqrt t)$, we have  $j_0,k-j_0 = \Theta (k)$.
    \item For $k = o(\sqrt t)$ we have  $k-j_0 = o(k)$.
    \item For $k = o(t^{1/4})$, we have  $j_0 = k$.
  \end{enumerate}
  Moreover, $|\W_i|\le (k+1)f(j_0)$ in general and $|\W_i|\le (1+o(1))f(k) = (1+o(1)){t\choose t-k}$ for $k = o(t^{1/4})$.
\end{observation}

We will need the following concrete statement concerning the size of $\T_n[\W_k].$

\begin{lemma}\label{peelbound}
  There is $\eps_0 > $ such that for any $\eps < \eps_0, n > n_0(\eps), t_0(\eps) < t \le n-2$ and $t^{0.01}\le k\le 2\epsilon t$ we have
  \begin{equation}\label{eqobs7}\frac{|\W_k|2^k}{n^k}\le 2^{-k}.\end{equation}
\end{lemma}
\begin{proof}
  Using the equality \eqref{eqfj} and the inequality ${m\choose \ell}\le (em/\ell)^\ell$, we get
  \begin{equation}\label{eqfj2} f(j)\le \Big(\frac{e^3tk^2}{j^3}\Big)^j(k+1)^{k-j}.\end{equation}
Let us first consider the case when $k\ge t^{0.99}$. Then we can apply Observation~\ref{w_bounds} (1) and get that the fraction in the inequality \eqref{eqfj2} for $f(j_0)$ is simplified to $((1+o(1))e^3k)^{j_0}$. We thus get
\begin{align*}  \frac{|\W_k|2^k}{n^k}\le& (k+1) \frac{f(j_0)2^k}{n^k} \\
\le&  \frac{((1+o(1))e^3k)^{j_0}(k+1)^{k-j_0+1}}{(n/2)^k}
\le  \frac{(100k)^k}{(n/2)^k} \le \ 2^{-k},\end{align*}
where the last inequality is due to $100k \le 200\eps t < 200\eps n \le n/4$ for $\eps$ small enough.

  Next, assume that $t^{0.01}\le k\le t^{0.99}$. In this regime, according to Observation~\ref{w_bounds}, we have $tk^2/j_0^3\le (1+o(1))t^{0.99}$, and thus we can conclude that \eqref{eqfj2} implies the following inequality.
  $$f(j)\le (Ct^{0.99})^k,$$
  where $C$ is some absolute constant. In this regime, we can conclude the following.
  \begin{align*}  \frac{|\W_k|2^k}{n^k}\le& (k+1) \frac{f(j_0)2^k}{n^k} \\
\le&    \frac{(k+1)(Ct^{0.99})^k}{(n/2)^k} \le \ 2^{-k},\end{align*}
where in the last inequality we use that for $n$ large enough $\frac{t^{0.99}}{n/2}\le 2/n^{0.01} < 1/100C$ and $(k+1)/100^k < 2^{-k}$.
\end{proof}

Finally, we estimate the same expression for small $k$.

\begin{lemma}\label{peelbound_small}
There is an absolute constant $C$ independent of $n$ such that for all $k < t^{0.01}$ we have $$\frac{|\W_k|2^k}{n^k}\le C.$$
\end{lemma}
\begin{proof}
    Observation \ref{w_bounds} (4) and the fact that $k = o(t^{1/4})$ imply $$(2/n)^k|\W_k| = O\left(\frac{2^k}{n^k}\cdot {t \choose t-k}\right) = O\left(\left(\frac{2et}{nk}\right)^k\right) = O(1).$$
\end{proof}

\section{Finding a dense piece}
We return to the notation used in Theorem~\ref{1SA}. The following theorem makes use of the peeling procedure and allows us to find a spanning forest $X$ such that the family $\ff(X)$ is much denser in $\T_n(X)$ than $\ff$ in $\T_n$.

\begin{theorem}\label{theorem_dense} There is such $\eps_0 > 0$ that for any $\eps < \eps_0, n > n_0(\eps)$ the following holds. Assume $t$ is a positive integer such that $2n^{1-\eps/3} \le t \le n-2, t > t_0(\eps)$.  Suppose that $\ff\subset \T_n$ is $t$-intersecting and satisfies $|\ff|\ge 2^{2-x}\cdot c_{n, t'}n^{n-2-t'}$ for some integer $x$, $t^{0.01}<x\le\frac 13\eps t$ , and where $t' = t- n^{1-\epsilon/3}$. Then there is a subset $X$ of size $t'$ such that $|\ff[X]|\ge \frac 12 {t'+x\choose t'}^{-1}|\ff|$.
\end{theorem}

\begin{proof}[Proof of Theorem~\ref{theorem_dense}]
Take a family $\ff$ as in the statement. First, we use Theorem~\ref{1SA} to get a spread approximation $\s$ of sets of size at most $q'=\min\{n-1, (1+\eps)t\}$, such that $\s$ is $t'$-intersecting, and a remainder $\ff'$, $|\ff'|\le n^{-\frac 13 \eps t}c_{n, t}n^{n-2-t}$, that satisfy $\ff\setminus \ff'\subset \T_n[\s]$. (It is easy to see that all requirements of Theorem~\ref{1SA} are met if we choose $\eps_0$ and $n_0(\eps)$ appropriately.) Note that $$|\ff'| \le n^{-\frac13\eps t}c_{n, t}n^{n - 2- t}  \overset{\ref{cnt_relative}, (i)}{\le} n^{-\frac13\eps t}c_{n, t'}n^{n - 2- t'}\cdot(2/n)^{t-t'} \le 2^{-\frac13 \eps t}c_{n,t'}n^{n - 2 - t'} \le $$ $$\le 2^{-x}c_{n, t'}n^{n-2-t'} \le |\ff|/4.$$ We determine two cases based on the value of $x$. 

\textbf{Case I: $x \le q' - t'$.} In this case we apply the peeling procedure from Section~\ref{sec_peel} to $\s$ with $t'$ playing the role of $t$, obtaining a sequence of families $\W_k$. We  want to bound the contribution of $|\ff[\W_k]|$ with large $k$ to $|\ff|$. To this end, we apply Lemma~\ref{peelbound} with $t'$ playing the role of $t$ for each $k= x+1,\ldots, q'-t'\le\epsilon t+(t-t')$. Note that all conditions are satisfied, most importantly, $k\le \epsilon t+(t-t')\le 2\epsilon t$. Then Lemma \ref{peelbound} implies the following bound on the contribution of layers $\W_k$ to $|\ff|$ for $k\ge x+1$.

\begin{align}\label{eq_sum} \sum_{k=x+1}^{q'-t'}|\ff[\W_k]| \le& \sum_{k=x+1}^{q'-t'}|\T_n[\W_k]| \le \sum_{k=x+1}^{q'-t'}|\W_{k}|\cdot c_{n, t'+k}n^{n-2-t'-k} \overset{\ref{cnt_relative}, (i)}{\le} \notag\\ \le& c_{n, t'}n^{n-2-t'}\sum_{k=x+1}^{q'-t'} \frac{|\W_k|2^k}{n^k}\overset{\ref{peelbound}}{\le} c_{n, t'}n^{n-2-t'}\sum_{k=x+1}^{q'-t'} 2^{-k}\le \notag\\\le& c_{n, t'}n^{n-2-t'}\cdot 2^{-x}. \end{align}

Put $\ff'' =(\ff\setminus \ff')\cap \T_n[\h_{t'+x}]$. (In words, we removed from $\ff$ the remainder $\ff'$ and `peeled off' the layers $\W_k$ of uniformity at least $t'+x+1$.) We have $|\ff|> 4c_{n, t'}n^{n-2-t'}\cdot 2^{-x},$ thus $|\ff\setminus\ff'|\ge \frac 34|\ff|$ and  $|\ff''|\ge \frac 12 |\ff|$. Take any set $F\in \h_{t'+x}$ (such set is defined correctly since $x + t' \le q'$). Via averaging, there is a $t'$-element subset $X\subset F$, such that $|\ff''(X)|\ge {t'+x\choose t'}^{-1}|\ff''|\ge \frac 12 {t'+x\choose t'}^{-1}|\ff|$.

\textbf{Case II: $x > q' - t'$.} Put $\ff'' = \ff \setminus \ff' \subset \ff[\s]$. Since $|\ff'| \le |\ff|/4,$ family $\s$ is not empty and $|\ff''| \ge |\ff|/2$. Take any set $S \in \s$, by averaging there is $X \subset S, |X| = t'$ such that $$|\ff''(S)| \ge {|S| \choose t'}^{-1}|\ff''| \ge \frac12{|S| \choose t'}^{-1}|\ff| \ge \frac12{t' + x' \choose t'}^{-1}|\ff''|.$$ The last inequality is implied by $|S| \le q' < x + t'$.

\end{proof}

The following corollary shall be the building block for a better spread approximation result.

\begin{corollary}\label{cor_dense}
There is such $\eps_0 > 0$ that for any $\eps < \eps_0, n > n_0(\eps)$ the following holds. Assume $t$ is a positive integer such that $2n^{1-\eps/7} \le t \le n-2, t > t_0(\eps)$.
Suppose that $\ff\subset \T_n$ is $t$-intersecting and satisfies $|\ff|\ge c_{n, t}n^{n-2-t} \cdot 2^{-\frac 12 n^{1-\epsilon/4}}$. Then there is a set $X$, $|X|\le t+n^{1-\epsilon/6}$, such that $\ff(X)$ is $10\epsilon t$-spread.
\end{corollary}

\begin{proof}
We will apply Theorem~\ref{theorem_dense} with $x = n^{1-\epsilon/4}$, choosing $\eps_0, n_0(\eps), t_0(\eps)$ properly. Since $$4c_{n, t'}n^{n - 2 - t'}\cdot 2^{-x} \le 4c_{n, t'}\cdot n^{n-2-t}\cdot n^{n^{1 - \eps/3}}\cdot2^{-x}  \overset{\ref{cnt_relative}, (iv)}{\le} 4c_{n, t}n^{n-2-t}\cdot n^{n^{1 - \eps/3}}\cdot n^{n^{1 - \eps/3}}\cdot2^{-x} \le$$ $$\le c_{n, t}n^{n-2-t}\cdot2^{2n^{1-\eps/3}\log_2n - n^{1 - \eps/4} + 2} \le c_{n, t}n^{n-2-t}\cdot 2^{-\frac12n^{1 - \eps/4}},$$ the condition on $|\ff|$ in Theorem~\ref{theorem_dense} is implied by $|\ff|\ge c_{n, t}n^{n-2-t}\cdot 2^{-\frac 12 n^{1-\epsilon/4}}$. It remains to check whether $t^{0.01} < x < \eps t/3$. It is true, because for $\eps_0$ small enough and $n_0(\eps)$ large enough$$t^{0.01} \le n^{0.01} \le n^{1- \eps/4} = x \le \frac{\eps}3n^{1 - \eps/7} \le \eps t/3.$$ Theorem~\ref{theorem_dense} gives us a $t'$-element set $X$, such that
$$|\ff(X)|\ge \frac 12 {t'+x\choose x}^{-1}|\ff|\ge n^{-x}|\ff|\ge c_{n, t'}n^{n-2-t'}\cdot 2^{-x}n^{-x} \ge$$$$\ge c_{n, t'}n^{n-2-t'}\cdot 2^{-x(\log_2n +1)} \ge c_{n, t'}n^{n-2-t'}\cdot 2^{-n^{1 - \eps/5}}.$$
We claim that there is a set $Y$ satisfying $|Y|\le t'+n^{1-\epsilon/6}$ and $X\subset Y$, such that, moreover, $\ff(Y)$ is $10\epsilon t$-spread. Indeed, find a maximum set $Y$, $Y\supset X$, that violates the $10\epsilon t$-spreadness of $\ff(X)$. If $|Y|> t'+n^{1-\epsilon/6}$, put $a = |Y| - t' > n^{1 - \eps/6}$.
$$|\ff(X)|\le (10\epsilon t)^{|Y|-t'}c_{n, |Y|}n^{n-2-|Y|} \le (10\epsilon t)^a\cdot 2^ac_{n, t'}n^{n-2-t'-a} \le  |\ff(X)|,$$
which is a contradiction. Here in the third inequality we used $$(20\eps t/n)^a \le 2^{-a}\le 2^{-n^{1-\eps/5}}$$ where the first inequality holds because we can choose $\eps < 1/40$ and the second one holds because $a > n^{1 - \eps/6}$.
\end{proof}

\section{Second spread approximation}

Using the corollary from the previous section, we can get a much better spread approximation of our family.

\begin{theorem}\label{2SA}
There is an $\eps_0 > 0$ such that for any $\eps < \eps_0, n > n_0(\eps)$ the following holds. Let $t\ge1$ be an integer such that $2n^{1-\eps/7} \le t \le n- 2, t > t_0(\eps)$. Consider a family $\ff\subset \T_n$ that is $t$-intersecting. Then there exists a  $t$-intersecting family $\mathcal S\subset \T_n^{\le t + n^{1 - \eps/6}}$ of spanning forests of size at most $t + n^{1 - \eps/6}$ and a family $\ff'\subset \ff$ such that the following holds.

  \begin{itemize}
    \item[(i)] $\ff\setminus \ff'\subset \T_n[\s]$;
    \item[(ii)] for any $B\in \s$ there is a family $\ff_B\subset \ff$ such that $\ff_B(B)$ is $10\epsilon t$-spread;
    \item[(iii)] $|\ff'|\le 2^{-\frac 12 n^{1-\epsilon/4}}c_{n, t}n^{n-2-t}$.
  \end{itemize}
\end{theorem}

\begin{proof}[Proof of Theorem~\ref{2SA}]
We iteratively apply Corollary~\ref{cor_dense}. We put $\ff_0 :=\ff$, and at the $i$-th step we get a set $Y$ and a family $\ff_i[Y]$. We then put $\ff_{i+1}:=\ff_i\setminus  \ff_i[Y]$ and apply Corollary~\ref{cor_dense} to $\ff_{i+1}$ as long as the corollary is applicable, that is, until step $N$ when $|\ff_N|\le c_{n, t'}n^{n-2-t'} \cdot 2^{-\frac 12 n^{1-\epsilon/4}}$ for some $N$. We then put $\ff':=\ff_N$. The collection of $Y$'s accumulated while running the procedure is our family $\mathcal S$, and we put $\ff_Y[Y]$ to be the corresponding family $\ff_i[Y]$.

The only property that we are left to verify is that $\s$ is $t$-intersecting. This is done in a way that is very similar to the proof of Lemma~\ref{lemtint_1SA}, and we sketch it below. Arguing indirectly, take $A_1,A_2\in \s$ that are not $t$-intersecting and, putting $x:=t-|A_1\cap A_2|$, for each $i\in [2]$ and $j=3-i$ define
 $$\g_i:=\ff_{A_i}(A_i)\setminus\Big\{F\in \ff_{A_i}(A_i): |F\cap A_j\setminus A_i|\ge x\Big\}.$$
 Using the bound on $|A_i|$ and the spreadness of $\ff_{A_i}(A_i)$, the size of the subtracted family is at most
 \begin{align*}{|A_j\setminus A_i|\choose x}\max_{X: |X| = x, X\cap A_i = \emptyset}|\ff_{A_i}(A_i\cup X)|\le&\\   {n^{1-\epsilon/6}+x\choose x}(10\epsilon t)^{-x}|\ff_{A_i}(A_i)|\le&\\
  (n^{1-\epsilon/6}+1)^x(10\epsilon t)^{-x}|\ff_{A_i}(A_i)|\le&\ \frac 12 |\ff_{A_i}(A_i)|.\end{align*}
  We conclude that $|\g_i|\ge \frac 12 |\ff_{A_i}(A_i)|$ and, consequently, that $\g_i$ is $5\epsilon t$-spread for both $i\in [2]$. Since $n$ is sufficiently large and $t>2n^{1-\epsilon/7}$, we have $5\epsilon t > 2^{5}\log_2(2n)$, and we can apply the coloring argument, finding two disjoint sets $U_1,U_2$, where $U_i\in \g_i$. Then $U_i\cup A_i$ violate the $t$-intersection property of $\ff$.
\end{proof}

\section{Proof of Theorem \ref{thmbigt}}

Fix such $\eps$ that it is smaller than both $\eps_0$ guaranteed by Theorems \ref{1SA} and \ref {2SA}. We will also assume $\eps < 0.01$. Let $n_0$ be a large number. We will determine how large it should be during the proof, at first we require it to satisfy the following conditions:
\begin{enumerate}
    \item $n_0 > 2^{19}$
    \item $n_0> n_0(\eps)$ where $n_0(\eps)$ also comes from Theorem \ref {2SA};
    \item for all $n \ge n_0$ we have $2n^{1 - \eps/7} < \frac n{4032\log_2n}$;
    \item for all $n \ge n_0$ Lemmas \ref{antr_relative}, \ref {ant1} hold.
\end{enumerate} Let $n, t$ be such integers that $n > n_0$ and $1 < t \le n-2$, let $\ff \subset \T_n$ be a $t$-intersecting family, $|\ff'| \ge c_{n, t}n^{n - 2 - t}$. We will prove that in fact equality holds and $\ff$ is the family of all spanning trees containing a fixed forest $F$ with $t$ edges and connected components of almost equal sizes.

If $t < 2n^{1-\eps/7}$ we are done by the result of \cite{EKR_trees}, so we may assume the opposite inequality holds. Apply Theorem \ref{2SA} and get a family $\s$ of sets of size at most $q:=t+n^{1-\eps/6}$ and a small remainder $\ff' \subset \ff, |\ff'| \le 2^{-\frac 12 n^{1-\epsilon/4}}|\ff| = o(|\ff|)$.

Next, apply the peeling procedure to $\s$. Fix $k$ to be the largest value such that, first, $k\le t^{0.01}$ and, second, $|\W_k|\ge t^{-0.5}(n/2)^k$. If there is no such $k$ put $k=0$. We stop the peeling procedure at $k$ and shall show that most of the family $\ff$ is contained in $\T_n[\W_k]$. Recall that $|\T_n(\W_k)|\le |\W_k|c_{n, t+k}n^{n-2-t-k}$.

First, we note that, by inequality \eqref{eqobs7}, as in the proof of Theorem~\ref{theorem_dense} in the series of inequalities \eqref{eq_sum}, we have \begin{equation}\label{eqres1}
\sum_{x = t^{0.01}+1}^{q-t}|\T_n(\W_x)|\le \sum_{x = t^{0.01}+1}^{q-t}2^{-x}c_{n, t}n^{n-2-t}\le 2^{-t^{0.01}}c_{n, t}n^{n-2-t}.                                        \end{equation}
Second, using the inequality on $|\W_x|$ for $x>k$ and (i) of Lemma~\ref{cnt_relative}, we have

\begin{equation}\label{eqres2}
\sum_{x = k+1}^{t^{0.01}}|\T_n(\W_x)|\le \sum_{x = k+1}^{t^{0.01}}|\W_x|c_{n, t+x}n^{n-2-t-x} \overset{\ref{cnt_relative}, (i)}{\le} c_{n, t}n^{n-2-t}\sum_{x = k+1}^{t^{0.01}}|\W_x|2^x/n^x \le  \end{equation}
$$\le c_{n, t}n^{n-2-t} \cdot t^{0.01} \cdot t^{-0.5} \le t^{-0.4}\cdot c_{n, t'}n^{n-2-t'}.$$

For the subsequent analysis, we need to better understand the structure of $\W_k$. Looking at the formula \eqref{eqfj} and the displayed equation above it, it should be clear that, (assuming $k<t^{0.01}$) for every $j\le k-1$ we have $f(j)=o\big(t^{-0.9}f(k)\big)= o\big(t^{-0.9}{t\choose t-k}\big)$. Moreover, if any two sets in $\W_k$ intersect in at least $t+1$ elements, then we can do the same analysis and calculations as the one that preceded \eqref{eqfj}, but with $A,B\in \W_k$ and with $I$ of size at least $t+1$ and $C$ intersecting $A,B$ in at least $t+1$ elements, and get a bound $|\W_k|=o\big(t^{-0.9}{t\choose t-k}\big)$. Thus, we conclude that, first, we have two sets $A,B\in \W_k$ such that $|A\cap B| = t$, and, moreover, by Observation~\ref{w_bounds}, most of the sets from $\W_k$ must intersect $U_0 = A\cap B$ in exactly $t-k$ elements and thus contain $A\Delta B$. 
  Let us fix one such set $C$ and put $I = A\cap B\cap C$. Let $I'$ stand for the set of elements that belong to exactly two out of $A,B,C$. That is, $A\cup B\cup C = I\sqcup I'$.

We are ready to complement inequalities \eqref{eqres1}, \eqref{eqres2} with the analysis of the sizes of layers below $k$. Put $\g_j:=\h_k\cap {{n \choose 2}\choose t+j}$ for $0\le j<k$.
First and foremost, let us bound the size of the intersection of $G\in \g_j$ with $I$. The set $G$ intersects each of $A,B,C$ in at least $t$ elements, thus $3t$ elements in total (with multiplicities). By double-counting, knowing that elements in $I$ have multiplicity $3$ and the elements in $I'$ have multiplicity $2$, we must have
$$|G|\ge |G\cap I|+|G\cap I'|\ge |G\cap I|+\Big\lceil\frac 32 (t-|G\cap I|)\Big\rceil= t+\Big\lceil\frac 12(t-|G\cap I|)\Big\rceil.$$ We have $|G|=t+j$, and  thus $t-|G\cap I|\le 2j$. In other words, $|G\cap I|\ge t-2j = |I|+k-2j$. From here, we can conclude that
$$|\g_j|\le \sum_{m=t-2j}^{t-k}{t-k\choose m}{3k\choose \big\lceil\frac 32 (t-m)\big\rceil}(k+1)^{j-\lceil (t-m)/2\rceil},$$
where the inequality is due to Lemma~\ref{lemma_peel}~(iii) that states that a subfamily of $\h_k(X)$ cannot be  $>k+1$-spread for a set $X$ of size at most $t+j-1$ (applied for $\g_j(X)$), and Observation~\ref{obs_spread_1}. Since $k\le t^{0.01}$, it is easy to see\footnote{Doing the same analysis as after \eqref{eqfj}} that the maximum in the sum above is attained when $m = t-2j$, and, moreover, the following holds
$$|\g_j|\le (1+o(1)){t-k\choose t-2j}{3k\choose 3j}(k+1)^{j-\lceil (t-(t-2j))/2\rceil}= (1+o(1)){t-k\choose t-2j}{3k\choose 3j}.$$

The number of sets on levels bellow $k$ is $$\sum_{j=0}^{k-1}|\T_n(\g_j)| \le (1+o(1))\sum_{j=0}^{k-1}{t-k \choose t-2j}{3k \choose 3j}c_{n, t+j}n^{n-2-t-j} \overset{\ref{cnt_relative}, (i)}{\le}$$ $$\le (1+o(1))c_{n, t}n^{n-2-t}\sum_{j=0}^{k-1}{t-k \choose t-2j}{3k \choose 3j}\cdot \frac{2^j}{n^j}.$$ If $t - 2j \ge t - k$ then $${t-k \choose t-2j}{3k \choose 3j}\cdot \frac{2^j}{n^j} = {3k \choose 3j}\cdot \frac{2^j}{n^j} \le (3k)^{3j}\cdot \frac{2^j}{n^j} \le \left(\frac{54k^3}{n}\right)^j \le \left(\frac{54}{n^{0.97}}\right)^j \le t^{-0.9}.$$ The opposite case is a little more complicated. Put $a = k - j \ge 1$, assume $j - a \ge 6$. $${t-k\choose t-2j}{3k\choose 3j}\cdot\frac{2^j}{n^j} = {t-k\choose j - a}{3k \choose 3a}\cdot\frac{2^j}{n^j} \le \left(\frac{et}{j - a}\right)^{j - a}(3^3k^3)^a\cdot\frac{2^j}{n^j}\le $$$$\le \left(\frac{2et}{(j-a)n}\right)^j \cdot \left(\frac{27k^3(j - a)}{et}\right)^{a}.$$ Since $j - a \ge 6 > 2e,$ then the first factor is less than $1$, thus $${t-k\choose t-2j}{3k\choose 3j}\cdot\frac{2^j}{n^j} \le \left(\frac{27k^3(j - a)}{et}\right)^{a} \le \frac{27k^4}{et} \le t^{-0.9}$$ where the last inequality holds for $t$ big enough. If $j - a \le 5$, we use a simpler bound on binomial coefficients. $${t-k\choose t-2j}{3k\choose 3j}\cdot\frac{2^j}{n^j} = {t-k\choose j - a}{3k \choose 3a}\cdot\frac{2^j}{n^j} \le t^{j - a}(3^3k^3)^a\cdot\frac{2^j}{n^j} \le$$ $$\le 2^j\cdot\left(\frac{27k^3}{t}\right)^a \le 32\cdot\left(\frac{54}{t}\right)^a \le t^{-0.9}.$$

From here we conclude that 
\begin{equation}\label{eqres3}
  |\T_n[\h_{k}\setminus \W_k]|\le \sum_{j = 0}^{k-1}|\T_n[\g_j]|\le (1+o(1))\cdot t^{0.01} \cdot t^{-0.9} c_{n, t}n^{n-2-t}\le t^{-0.7}c_{n, t}n^{n-2-t}.
\end{equation}

From inequalities \eqref{eqres1},\eqref{eqres2} and \eqref{eqres3}, together with the bound on $|\ff'|$ that is given by the application of Theorem~\ref{2SA}, we conclude that, as long as $|\ff| \ge c_{n, t}n^{n-2-t}$, all but an $o(1)$-fraction of all sets in $\ff$ are contained in $\T_n(\W_k)$. Using the last conclusion of Observation~\ref{w_bounds}, all but a $o(1)$-fraction of the sets in $\W_k$,  are contained in the union $A\cup B$, where $A,B\in \W_k$ and $|A\cap B| = t$. (cf. the paragraph after \eqref{eqres2}) That is, they belong to ${A\cup B\choose t+k}$, where $|A\cup B| =t+2k$.

Put $$\ff_1 = \bigcup_{C\subset A\cup B, |C|=t+k}{\ff[C]} \text{ and } \ff_2 = \ff \setminus \ff_1.$$

The family $\ff_2$ consists of some sets from $\ff \setminus \ff[\W_k]$ and some sets from $\cup_{C\in \W_k, C\not\subset A \cup B} \ff[C]$. There are no more than $o(|\ff|)$ sets of the first type. The number of sets of the second type does not exceed $$o\left(|\W_k|c_{n, t+k}n^{n-2-t-k}\right) = o\left(2^k|\W_k|c_{n, t}n^{n-2-t-k}\right) =$$$$= o\left(\aaa_{n, t, 0}\cdot\left(\frac2n\right)^k\cdot|\W_k|\right) = o(|\ff|) \cdot \left(\frac2n\right)^k\cdot|\W_k|,$$ where the first inequality is due to (i) of Lemma \ref{cnt_relative} and the third one is due to the fact that $\ff$ is the largest $t$-intersecting family.

Finally, Lemma \ref{peelbound_small} gives $(2/n)^k|\W_k| = O(1)$. We  conclude that all but $o(|\ff|)$ sets in the family $\ff$ are contained in $\uu_{n, t, k, A}$. Note that the rate of decay of $o(\cdot)$ is polynomial in  $n$. Lemma \ref{antr_relative} implies that we must have $k = 0$ or $k = 1$ due to extremality of $\ff$. Let us  summarize what we have achieved  by now.

\begin{theorem}\label{bigt_asymp}
    For all $n \ge n_0$ and $2 \le t \le n-2$ the following holds. If $\ff\subset \T_n$ is $t$-intersecting then $|\ff| \le (1 + o(1))c_{n, t}n^{n-2-t}$, where $c_{n, t}$ is the largest possible product of $n-t$ positive integers with sum $n$. Moreover, there is a set of edges $A$ such that all but an $o(1)$-fraction of sets in $\ff$ lie in $\uu_{n, t, r, A}$, where $r \in \{0,1\}$.
\end{theorem}
Theorem \ref{bigt_asymp} gives us an 'asymptotic' version of the result we aim to proof. We are left to show that, in order to be extremal, the remainder, not lying in the corresponding $\uu$, is empty. We treat several cases depending on $k$ and the size of the remainder.

\textbf{Case Ia: $k = 0, |\ff_2| \ge {n^{-\eps t/4}}c_{n, t}n^{n-2-t}$.} In this case $|A| = t, |\ff_1| = (1+o(1))|\ff| = (1+o(1))c_{n, t}n^{n-t-2}$. According to the Theorem \ref{1SA} there is a forest $S$ of size not exceeding $(1+\eps)t$ and a family $\ff_S\subset \ff_2$ such that $\ff_S(S)$ is $n^{\eps/2}$-spread. We start by deleting from $\ff_1$ all sets that have at least $t$ common elements with $S$. Denote $\ff_3$ the family obtained from $\ff_1$ during this procedure.

First note that $A\cap S \neq A$, otherwise $F_S[S] \subset \ff_1$. Put $z = |A\cap S| < t.$ We should delete no more than $x = {(1+\eps)t - z \choose t - z}c_{n, 2t-z}n^{n - 2 - 2t + z}$ sets since each of the sets that we are about to delete contains all $t$ elements of $A$ and at least some $t-z$ elements from $S\setminus A$. First consider $t - z \ge 2$. If $t \ge 2n/3$ a simple bound on the binomial coefficients together with (ii) of Lemma \ref{cnt_relative} gives $$x \le \frac{((1+\eps)n)^{t-z}}{(t-z)!}\cdot n^{n-2-2t+z} \cdot c_{n, t+t-z} \leqslant c_{n, t}n^{n-2-t}\frac{((1+\eps)e/3)^{t-z}}{(t-z)!} \le (1+o(1))|\ff|/2,$$ where the last inequality holds for $\eps < 0.1$. If $t < 2n/3$ then $(1+\eps)t - z \leqslant (1 + \eps)2n/3$ and using the same bound and (i) of Lemma \ref{cnt_relative} we get $$x \leqslant c_{n, t}n^{n-2-t}\frac{((1+\eps)4/3)^{t-z}}{(t-z)!} \le (1+o(1))\cdot 0.9|\ff|,$$ where the last inequality holds for $\eps < 0.01$.

In the case $t-z = 1$ the bound on $x$ can be rewritten as follows: $$x \le (\eps t+1)\cdot c_{n, t+1}\cdot n^{n-2-t-1}\le (1 + o(1)) \frac{2\eps t}{n}\cdot |\ff| \leqslant (1+o(1))|\ff|/2,$$ where in the second inequality we used (i) of Lemma \ref{cnt_relative} and the third holds for $\eps < 1/4$. Since we could have chosen $n_0$ properly, $|\ff_3| \ge (1+o(1))|\ff|/10 \ge |\ff|/11$. 

Since $\ff$ was maximal and $|\ff_3|$ is only smaller than $|\ff|$ by a constant factor, $\ff_3$ is big enough and we can once again use the Theorem \ref{1SA}. We obtain such $T$ and $\ff_T\subset \ff_1$ that $\ff_T(T)$ is $n^{\eps/2}$ spread. Without loss of generality we may assume $A \subset T$. 

The procedure we used to construct $\ff_3$ implies $|S\cap T| < t$. It remains to use Lemma \ref{find_disjoint} with $R = n^{\eps/2}$ and families $\ff_S(S)$ and $\ff_T(T)$ of uniformity less than $n$. Since $R = n^{\eps/2} > 2^{5}\log_2(n-1)$ for $n$ large enough, there are pairwise disjoint $G_T \in \ff_T(T), G_S \in \ff_S(S)$. Trees $A_1 = G_T \cup T$ and $A_2 = G_S\cup S$ both lie in $\ff$. Thus they must $t$-intersect. However, $$|A_1 \cap A_2|  = |A_1 \cap S| + |A_1 \cap G_S| = |A_1 \cap S| < t, $$ where the last inequality holds because $A_1 \in \ff_3$ and $\ff_3$ does not contain sets having $t$ or more common elements with $S$. This leads us to contradiction.

\textbf{Case Ib: $k = 0, |\ff_2| < {n^{-\eps t/4}}c_{n, t}n^{n-2-t}$.} We will use Lemma \ref{avoiding}, however it is only applicable if $t \le n - 55$. However it is not an important restriction since for large $n$ and $n - 55 < t \le n-2$ the following holds: $$|\ff_2| < {n^{-\eps t/4}}c_{n, t}n^{n-2-t} \le n^{-\eps t/4} \cdot n^{n - t} \cdot n^{n - 2 - t} \le n^{-\eps t/4 + 110 - 2} = o(1).$$ The second inequality is true because $c_{n, t}$ is the maximal possible product of $n-t$ numbers with sum $n$ which can bu roughly bounded by $n^{n-t}$. Thus we may assume $t \le n - 55$. 

In this case for $n$ big enough holds $$|\ff_2| < \left(\frac{2}{n^{\eps/4}}\right)^tn^{n-2-t} < n^{n-t-27}.$$
If $\ff_2$ contains at least one tree that is not a star, then Lemma \ref{avoiding} implies $$|\ff_1|\leqslant c_{n, t}n^{n-2-t} - n^{n-t-27}.$$ Thus $|\ff| < \aaa_{n, t, 0}$ which contradicts our assumption that $\ff$ is maximal.

It remains to treat the case when $\ff_2$ is not empty and consists only of stars. Here $|\ff_2|\le n$. Take one of these stars $S$. However due to Lemma \ref{avoiding_star} and the fact that $t \le n - 55$ there are more than $n$ trees that contain $A$ and has less than $t$ common edges with $S$. Thus any case where $k=0$ and $\ff_2$ is not empty leads to contradiction. The only remaining case is when $\ff = \T_n[F]$ where $F$ is a forest consisting of $n-t$ connected components of almost equal sizes.

\textbf{Case IIa: $k = 1, |\ff_2| \ge {n^{-\eps t/4}}c_{n, t}n^{n-2-t}$.} Lemma \ref{antr_relative} implies that this case is only relevant if $t < n/2$. Once again we start with finding a forest $S$ and a family $\ff_S \subset \ff_2$ such that $\ff_S(S)$ is $n^{\eps/2}$ spread and $|S| \leq (1+\eps)t$. Then we construct $\ff_3$ from $\ff_1$ by eliminating all trees that share at least $t$ edges with $S$.

Recall that there is such a set $A$ of size $t+2$ that $\ff_1 \subset \cup_{B\subset A, |B| = t+1}\T_n[B]$. Put $z = |A \cap S| < t$. We should delete
no more than $$x = (t+2){(1+\eps)t-z \choose t-z}\cdot c_{n, t+1+t-z}n^{n-2-(t+1+t-z)}$$ sets. The expression above arises because we should choose $B$ --- one of $t+2$ subsets of size $t+1$ in $A$, the choose $t-z$ elements from $S\setminus A$ and finally choose a tree that contains $B$ and those $t-z$ vertices from $S\setminus A$. Once again the analysis depends on $t-z$ and $t$. First we consider $t - z = 1$. Then $x$ can be rewritten as 
$$x = (t+2)(1 + \eps t)\cdot c_{n, t+2}n^{n-2-t - 2} \le (1+o(1))c_{n, t}n^{n-2-t}\frac{4\eps t^2}{n^2} \le $$$$ \le (1+o(1))|\ff| \cdot \eps \le (1+o(1))\cdot |\ff|/2.$$
Second we consider $t-z \ge 2$. Due to a bound on binomial coefficients and (i) of Lemma \ref{cnt_relative} we get $$x \le (t+2)\frac{((1+\eps)/2)^{t-z}n^{t-z}}{(t-z)!} \cdot 2^{t+1-z}c_{n, t}n^{n- 2 - 2t + z - 1} \le $$ $$\le c_{n, t}n^{n-2-t}\frac{2(t+2)(1+\eps)^{t-z}}{n(t-z)!} \le (1+o(1))\cdot 3|\ff|/4$$

In both cases $|\ff_3| \ge (1+o(1))|\ff|/4$ and we finish the proof in the exactly same way as in the case Ia, since once again $|\ff_3|$ differs from $|\ff|$ only by a constant factor.

\textbf{Case IIb: $k = 1, |\ff_2| < {n^{-\eps t/4}}c_{n, t}n^{n-2-t}$.} By analogy with the case Ib we get $|\ff_2| < n^{n-(t + 1)-27}$.  Note that in this case we do not need to state the bound $t \le n - 55$ explicitly, because it can be guaranteed by a proper choice of $n_0$ since $t< n/2$. Fix arbitrary $B\subset A, |B| = t+1$. If $\ff_2$ contains at least one tree $T$ that is not a star, then $|\ff[B]| \leqslant |\ff| - n^{n - (t+1)-27}$ which leads to contradiction. If $\ff_2$ is not empty but consists only of stars, then $|\ff_2| \leqslant n$ and we come to a contradiction using Lemma \ref{avoiding_star}.

Therefore $\ff \subset  \cup_{B\subset A, |B| = t+1}\T_n[B]$. Since the RHS is $t$-intersecting, we get $\ff = \uu_{n, t, 1, A}$. This family can not be extremal provided Lemma \ref{ant1} is true. 

Thus the only optimal example is trivial, which concludes the proof of the main theorem.

\section{Acknowledgements} This work was in part done during LIPS'2025, a summer research program for students at MIPT.


\begin{thebibliography}{100}
\bibitem{AK} R. Ahlswede and L.H. Khachatrian, {\it The Complete Intersection Theorem for Systems of Finite Sets}, European Journal of Combinatorics. 18 (1997), 125--136.


\bibitem{Alw} R. Alweiss, S. Lovett, K. Wu, and J. Zhang, {\it Improved bounds for the sunflower lemma}, arXiv:1908.08483, (2019)

\bibitem{EFP} D. Ellis, E. Friedgut, and H. Pilpel, {\it Intersecting Families of Permutations}, J. Amer. Math. Soc. 24 (2011), 649--682.

\bibitem{EKR} P. Erd\H os, C. Ko, and R. Rado, \textit{Intersection theorems for systems of finite sets}, The Quart. J. Math. 12 (1961), N1, 313--320.

\bibitem{Film} Y. Filmus, {\it The weighted complete intersection theorem}, Journal of Combinatorial Theory, Series A 151 (2017), 84--101. 



\bibitem{F1} P. Frankl, \textit{The Erd\H os-Ko-Rado theorem is true for n=ckt}, Combinatorics (Proc. Fifth Hungarian Colloq., Keszthely, 1976), Vol. I, 365--375, Colloq. Math. Soc. J\'anos Bolyai, 18, North-Holland.


\bibitem{EKR_trees} P. Frankl, G. Hurlbert, F. Ihringer, A. Kupavskii, N. Lindzey, K. Meagher, V. R. T. Pantangi, {\it Intersecting Families of Spanning Trees, arXiv preprint  arXiv:2502.08128}, (2024).


\bibitem{FF3} P. Frankl and Z. F\"uredi, \textit{Beyond the Erdos-Ko-Rado theorem},  J. Combin.  Theory Ser. A  56 (1991) N2, 182--194.

\bibitem{FKhr} P. Frankl, A. Kupavskii, {\it The Hajnal--Rothschild problem}, arXiv:2502.06699 (2024)
\bibitem{FW86} P. Frankl, R.M. Wilson, {\it The Erd\H os-Ko-Rado theorem for vector spaces}, Journal of Combinatorial Theory, Series A 43 (1986), 228--236. 


\bibitem{Hu} L. Hu, {\it Entropy Estimation via Two Chains: Streamlining the Proof of the Sun-
flower Lemma} \url {https://theorydish.blog/2021/05/19/entropy-estimation-via-two-chains-streamlining-the-proof-of-the-sunflower-lemma/}, (2021).

\bibitem{Keevash} P. Keevash, N. Lifshitz, E. Long, D. Minzer, {\it Global hypercontractivity and its applications}, arXiv preprint  arXiv:2103.04604, (2021)

\bibitem{Kat} G. Katona, {\it Intersection theorems for systems of finite sets}, Acta Math. Acad. Sci. Hungar. 15 (1964) 329--337.

\bibitem{KLMS} N. Keller, N. Lifshitz, D. Minzer, and O. Sheinfeld, {\it On $t$-intersecting families of permutations}, http://arxiv.org/abs/2303.15755v2

\bibitem{Kup55} A. Kupavskii, {\it Intersection theorems for uniform subfamilies of hereditary families} (2023), arXiv.2311.02246.


\bibitem{permutations} A. Kupavskii, {\it An almost complete t-intersection theorem for permutations, arXiv preprint arXiv:2405.07843}, (2024).


\bibitem{Fedya} A. Kupavskii, F. Noskov, {\it Linear dependencies, polynomial factors in the Duke--Erd\H os forbidden sunflower problem}, arXiv preprint arXiv:2410.06156, (2024)


\bibitem{spreads} A. Kupavskii and D. Zakharov, {\it Spread approximations for forbidden intersections problems}, Advances in Mathematics, 445:109653, (2024).

\bibitem{tree_counter} Linyuan Lu, Austin Mohr, and L\'aszl\'o  Sz\'ekely, {\it Quest for negative dependency graphs}, Recent advances in harmonic analysis and applications, volume 25
of Springer Proc. Math. Stat., pages 243–258. Springer, New York, (2013).


\bibitem{Sto} M. Stoeckl, {\it Lecture notes on recent improvements for the sunflower lemma}, \url {https://mstoeckl.com/notes/research/sunflower_notes.html}

\bibitem{W} R.~M. Wilson, {\it The exact bound in the Erd\H s--Ko--Rado theorem}, Combinatorica 4 (1984), N2-3, 247--257.

\end{thebibliography}
\end{document}